\documentclass[11pt,a4paper]{amsart}
\usepackage{color}
\usepackage{amssymb}
\def\R{\mathbb{R}}
\def\erre{\mathbb{R}}
\def\m1{{I\!\!M}}
\newcommand{\rdue}{\erre^2}

\def\sideremark#1{\ifvmode\leavevmode\fi\vadjust{\vbox to0pt{\vss% the remark
 \hbox to 0pt{\hskip\hsize\hskip1em%                          will appear only
 \vbox{\hsize2.1cm\tiny\raggedright\pretolerance10000%          on the side
  \noindent #1\hfill}\hss}\vbox to15pt{\vfil}\vss}}}%

\newcommand{\de}{\omega}
\newcommand{\grad}{\nabla}

\renewcommand{\to}{\rightarrow}
\newcommand{\pa}{\partial}

\newcommand{\ino}{\int_{\Omega}}

\newcommand{\ainf}{\mbox{as\;}\;n\to+\infty}

\newcommand{\fo}{\forall\,}

\newcommand{\rife}[1]{(\ref{#1})}
\newcommand{\ov}[1]{\overline{#1}}

\newcommand{\sscp}{\scriptscriptstyle}
\newcommand{\dsp}{\displaystyle}

\renewcommand{\dfrac}{\displaystyle\frac}
\newcommand{\fineproof}{\hspace{\fill}$\square$}

\renewcommand{\i}{\infty}

\renewcommand{\a}{\epsilon}

\newcommand{\eps}{\varepsilon}
\newcommand{\dt}{\delta}

\newcommand{\al}{\alpha}

\newcommand{\sg}{\sigma}
\newcommand{\ga}{\gamma}
\newcommand{\om}{\Omega}
\newcommand{\lm}{\lambda}
\newcommand{\vp}{\varphi}

\newtheorem{theorem}{Theorem}[section]
\newtheorem{proposition}[theorem]{Proposition}
\newtheorem{lemma}[theorem]{Lemma}
\newtheorem{corollary}[theorem]{Corollary}
\newtheorem{remark}[theorem]{Remark}
\newtheorem{definition}[theorem]{Definition}
\newcommand{\brm}{\begin{remark}}
\newcommand{\erm}{\end{remark}}
\newcommand{\bdf}{\begin{definition}}
\newcommand{\edf}{\end{definition}}
\newcommand{\bte}{\begin{theorem}}
\newcommand{\ete}{\end{theorem}}
\newcommand{\bpr}{\begin{proposition}}
\newcommand{\epr}{\end{proposition}}
\newcommand{\ble}{\begin{lemma}}
\newcommand{\ele}{\end{lemma}}
\newcommand{\bco}{\begin{corollary}}
\newcommand{\eco}{\end{corollary}}
\newcommand{\beq}{\begin{equation}}
\newcommand{\eeq}{\end{equation}}
\newcommand{\bdm}{\begin{displaymath}}
\newcommand{\edm}{\end{displaymath}}

\newcommand{\graf}[1]{\left\{\begin{array}{ll}#1\end{array}\right.}

\newcommand{\bel}{\begin{equation}}
\newcommand{\eel}{\end{equation}}

\begin{document}
\numberwithin{equation}{section}
\parindent=0pt
\hfuzz=2pt
\frenchspacing

\title[Global existence for a supercritical Keller-Segel type system]{A global existence result for a Keller-Segel type system with supercritical initial data}

\author[D.B. \& D.C.]{Daniele Bartolucci$^{(1,\dag)}$, Daniele Castorina$^{(2,\ddag)}$}

\thanks{2010 \textit{Mathematics Subject classification:} 35J61, 35K45, 35K57, 35K58}

\thanks{$^{(1)}$Daniele Bartolucci, Department of Mathematics, University
of Rome {\it "Tor Vergata"}, \\  Via della ricerca scientifica n.1, 00133 Roma,
Italy. e-mail:bartoluc@mat.uniroma2.it}

\thanks{$^{(2)}$Daniele Castorina, Dipartimento di Matematica, Universit\` a di Padova, \\
Via Trieste 63, 35121 Padova, Italy. e-mail: castorin@math.unipd.it}

\thanks{$^{(\dag)}$Research partially supported by FIRB project {\sl
Analysis and Beyond} and by MIUR project {\sl Metodi variazionali e PDE non lineari}} 

\thanks{$^{(\ddag)}$Research partially supported by project {\sl Bando Giovani Studiosi 2013 - Universit\` a di Padova - GRIC131695}}

\begin{abstract}
We consider a parabolic-elliptic Keller-Segel type system, which is related to a simplified model of 
chemotaxis. Concerning the maximal range of existence of solutions, there are essentially two kinds of results: 
either global existence in time for general subcritical ($\|\rho_0\|_1<8\pi$) initial data, or blow--up in finite time for suitably chosen 
supercritical ($\|\rho_0\|_1>8\pi$) initial data 
with concentration around finitely many points. As a matter of fact there are no results claiming the existence of global solutions in the 
supercritical case. We solve this problem here and prove that, for a particular set of initial data which share large supercritical masses, 
the corresponding solution is global and uniformly bounded.
\end{abstract}
\maketitle
{\bf Keywords}: Keller-Segel; Chemotaxis; Global solution; Supercritical problems.

\section{Introduction and main result}\label{intro}

Let $\om\subset \R^2$ be any smooth and bounded domain, we  consider the following parabolic-elliptic Keller-Segel type system 
\begin{equation}\label{system}
\begin{cases}
\rho_t = \nabla \cdot (\nabla \rho - \rho \nabla (u + \log V)), \quad &x \in \Omega, t >0\\
-\Delta u = \rho, \quad &x \in \Omega, t \geq 0\\
\rho (x,0) = \rho_0 (x)\geq 0 , \quad \int_{\Omega} \rho_0 = \lambda \quad &x \in \Omega\\
\frac{\partial \rho}{\partial \nu} - \rho \frac{\partial (u + \log V)}{\partial \nu} = 0, \quad u = 0 \quad &x \in \partial \Omega, t >0\\
\end{cases}
\qquad P(\lambda,\Omega)
\end{equation}

We assume that $V$ satisfies,
\begin{equation}\label{hypV}
V \in C^{0,1} (\Omega) \, \text{ and } \, 0 < a \leq V(x) \leq b.
\end{equation}

This is a simplified version of a chemotaxis model first introduced in \cite{KS}. The analysis of these kind of problems
has attracted a lot of attention in recent years and we refer 
the interested reader to the monograph \cite{suzC} for a complete account about this topic.\\

We say that $(\rho,u)$ is a classical solution of $P(\lambda,\Omega)$ in $[0,T]$ if $\rho\geq 0$,
\beq\label{global}
\rho \in C^{0}(\,\ov{\om}\times[0,T] )\cap C^{2,1}\left(\, \ov{\om}\times (0,T] \,\right),\quad u\in C^{2,1}\left( \, \ov{\om}\times (0,T] \,\right),
\eeq
and $(\rho,u)$ solves \rife{system}. 
We say that $(\rho,u)$ is a global solution of $P(\lambda,\Omega)$ if it is a classical solution in $[0,T]$ for any $T>0$.\\

Concerning the maximal range (in time) of existence of solutions to $P(\lambda,\Omega)$,  
there are essentially two kinds of results. The first one yields sufficient conditions 
to guarantee that the solution is global and uniformly bounded, that is
\beq\label{unifbd}
\sup\limits_{t>0}\sup\limits_{x\in\ov{\om}}\rho(x,t)\leq C.
\eeq
In fact it is well known that if $\rho_0$ is smooth and $\lm<8\pi$ then $P(\lambda,\Omega)$ admits a unique solution which is 
global and uniformly bounded, see \cite{suz} and also \cite{Bil}, \cite{GZ}. The second class of results is about sufficient conditions 
which guarantee that blow up occurs in finite time, 
that is, there exist $T_{\rm max}>0$ such that
\beq\label{unifbd1}
\lim\limits_{t\nearrow T^{-}_{\rm max}}\sup\limits_{x\in\ov{\om}}\rho(x,t)= +\infty.
\eeq
These conditions require the initial density $\rho_0$ to satisfy $\lm>8\pi$ and to be "peaked" around some point,  see 
\cite{SeSu}, \cite{suz} and also \cite{JL}.
Other intermediate situations may occur, such as for example blow up in infinite time when $\lm=8\pi$, see \cite{suz}.
So the value $\lm=8\pi$ is said to be the "critical" threshold and the study of \rife{system} is generally divided in the subcritical and supercritical 
regime according to whether $\lm<8\pi$ or  $\lm>8\pi$ respectively.
Of course, there are many other results which are concerned with the blow up rate and the structure of the blow up set, see \cite{suzC} and more recently \cite{suz} 
for further details.\\
As a matter of fact there are no results at hand claiming the existence of global solutions in the supercritical case $\lm>8\pi$. 
We solve this problem here and prove that \rife{system} may admit global and uniformly bounded solutions in the supercritical case as well. 
In fact we are able to find a particular set of initial data which share arbitrarily large supercritical masses $\lm>8\pi$ such that the corresponding solution is global and uniformly bounded. More exactly we have the following  global existence result for \rife{system}.\\

\begin{theorem}\label{kstype}

{\rm (a)} Let $\om$ be any smooth, bounded and simply connected domain.
For any $c\in(0,1]$ let $D = \frac{a}{b}$ (as given in \eqref{hypV}) and $c_{\sscp D} = c D$. Then there exist $\overline{\a}_*>\underline{\a}_*( c_{\sscp D} )>0$ 
such that
if $\{\a^2x^2+y^2\leq \beta_-^2\}\subset\Omega\subset\{\a^2x^2+y^2\leq\beta_+^2\}$ with $c=\tfrac{\beta_-^2}{\beta_+^2}$ then, for
any $\a\in(0,\underline{\a}_*(c_{\sscp D})]$ and for any $\lm\leq\lm_{\a,c_{\sscp D}}$, there exist initial data $\rho_{0}$ such that 
$P(\lambda,\Omega)$ admits a unique 
global and uniformly bounded solution $(\rho_\lambda, u_\lambda)$. Here $\underline{\lm}_{\a,c_{\sscp D}}< \lm_{\a,c_{\sscp D}}<\overline{\lm}_\a$ and 
$\underline{\lm}_{\a,c_{\sscp D}}$, $\overline{\lm}_\a$ are strictly decreasing (as functions of $\a$) in $(0,\underline{\a}_* (c_{\sscp D})]$, $(0,\overline{\a}_*]$ 
respectively with $\underline{\lm}_{\a,c_{\sscp D}}\simeq \frac{4\pi c_{\sscp D}}{(8-c_{\sscp D})\a}$, $\ov{\lm}_{\a}\simeq \frac{2\pi}{3 \a}$ as $\a\to 0^+$.\\

{\rm (b)}
There exists $\bar N> 4\pi$ such that if $\Omega$ is any open, bounded and convex set 
whose isoperimetric ratio, $N\equiv N(\Omega)=\frac{L^2(\pa\Omega)}{A(\Omega)}$, satisfies $N\geq\bar N$, then
for any $\lm\leq\lm_\textnormal{\tiny{$N$}}$ there exist initial data $\rho_{0}$ such that $P(\lambda,\Omega)$ admits a unique 
global and uniformly bounded solution $(\rho_\lambda, u_\lambda)$. Here
$\underline{\Lambda}_\textnormal{\tiny{$N$},D}<\lambda_\textnormal{\tiny{$N$}}<
\overline{\Lambda}_\textnormal{\tiny{$N$}}$ and $\underline{\Lambda}_\textnormal{\tiny{$N$},D}$, $\overline{\Lambda}_\textnormal{\tiny{$N$}}$
are strictly increasing in $N$ and $\underline{\Lambda}_\textnormal{\tiny{$N$},D}\simeq\frac{D \pi^2 N}{16(32-D)}+\mbox{\rm O}(1)$,
$\overline{\Lambda}_\textnormal{\tiny{$N$}}\simeq\frac{2\sqrt3 N}{\pi}+\mbox{\rm O}(1)$ as $N\to+\infty$.
\end{theorem}

The proof of Theorem \ref{kstype} is based on the following observation. System \eqref{system} admits a natural Lyapunov functional, which is 
the free energy \rife{free} below.  Suppose that we were able to find a strict local free energy minimizer, say $\rho_{0,\lm}$.
Then, in a carefully defined dual topology, those solutions of \eqref{system} with initial data in a small enough neighbourhood of $\rho_{0,\lm}$ should be trapped 
there for any $t>0$.

This information should yield the uniform estimates needed to prove global existence as well as uniform bounds. So the problem is to find out such minimizers 
and a good topology to work with.  
We seek these kind of minimizers in the class 
of stationary states $(\rho,u)$ of \eqref{system} which therefore satisfy
\begin{equation}\label{steady}
\begin{cases}
-\Delta \rho = \nabla \cdot (\rho \nabla (u + \log V)), \quad &x \in \Omega\\
-\Delta u = \rho, \quad &x \in \Omega\\
\int_{\Omega} \rho = \lambda \\
\frac{\partial \rho}{\partial \nu} - \rho \frac{\partial (u + \log V)}{\partial \nu} = 0, \quad u = 0 \quad &x \in \partial \Omega\\
\end{cases}
\end{equation}

We choose a particular steady state $(\rho_{0,\lm}, u_\lambda)$ of the 
form $\rho_{0,\lm} = \lambda \frac{V e^{u_\lambda}}{\int_\Omega V e^{u_\lambda}}$ (which satisfies the Neumann type 
boundary condition in \eqref{steady} automatically) and so reduce the problem to the existence for large $\lm$ of a free energy minimizer in 
the form of a (possibly weak) solution 
$u_\lambda$ of the following mean field equation with homogeneous Dirichlet boundary conditions:
\begin{equation}\label{mf}
\begin{cases}
-\Delta u_\lambda = \lambda \dfrac{V e^{\dsp u_\lambda}}{\int_\Omega V e^{\dsp u_\lambda}} \quad\mbox{in }\Omega\\
\quad u_\lambda = 0 \quad \mbox{on } \partial \Omega
\end{cases}
\end{equation}

The existence of free energy minimizers taking the form $\lambda \frac{V e^{u_\lambda}}{\int_\Omega V e^{u_\lambda}}$ with $\lm<8\pi$ is well 
known \cite{clmp2} where $u_\lm$ is a minimizer of the corresponding variational functional, see \rife{energy} below. 
On the contrary, in case $\lm>8\pi$ solutions of \rife{mf} on domains with non trivial topology are well known to exist \cite{CLin2} which are 
not minimizers of \rife{energy} in general. The reason behind this issue is that if $\lm>8\pi$ then both the variational functional \rife{energy} and the 
free energy \rife{free} are not bounded from below. Therefore, in particular, any local minimizer for $\lm>8\pi$ won't be a global one and would correspond not 
to a stable state but in fact to a so called metastable state from the dynamical point of view. 
Luckily enough, the existence for large $\lm$ of solutions of \rife{mf} on narrow domains 
which minimize the functional \rife{energy} has been recently established in  \cite{BdM2} in case $V$ is constant and then generalized to 
the case of non constant $V$ for Liouville systems in \cite{BCV}. In \cite{BCV} it has also been shown that these solutions naturally yield 
minimizers of the multidimensional analogue of the free energy \rife{free} in a suitable dual Orlicz-type topology. However we face a more subtle problem here 
since some properties which are almost obvious in the $H^1_0(\om)$-topology (see  \rife{220415.1}) become more delicate in the Orlicz setting.
As a consequence some care is needed to show that local minimizers of the free energy inherits that property from local minimizers of \rife{energy} 
in a suitable large enough sub cone, see in particular {\rife{200415.1}} in Proposition \ref{min} below. 
To make the exposition self contained we will also provide the details of the existence result, see Theorem \ref{subsup}.
This might be also useful since it clarifies the role played by a non constant $V$ as well as the minimizing properties of those solutions in the scalar case.\\
Once we have found the desired minimizers with the necessary topological informations, 
then we can prove that in fact some crucial a priori bounds hold, see \rife{main1} below. By using these a priori 
estimates, then the proof of the global existence result could follow in principle from an adaptation of 
well known arguments \cite{Bil}, \cite{GZ} to the analysis of system \rife{system}. It turns out that we do not need to work out this 
argument since this kind of adaptation has been recently worked out in \cite{suz} via a beautiful shorter proof of uniform boundedness and global 
existence of solutions based on \rife{main1}.   

\brm
It is likely that, by using a uniqueness result obtained in \cite{BdM2}, one could also prove that the $\rho_\lm$ found in Theorem \ref{kstype} 
converge as $t\to+\infty$ to the pair $(\rho_{0,\lm},u_{\lm})$ solving \rife{mf}.
\erm

This paper is organized as follows. In section  \ref{s1} we recall some well known facts about Keller-Segel systems. 
In section  \ref{s2} we prove a global existence result for \eqref{system}. In section \ref{s3} we obtain the needed free energy minimizers with large masses.
Finally, in section \ref{theoproof}, we prove Theorem \ref{kstype}. 

\bigskip
\bigskip

\section{Preliminaries}\label{s1}
We collect here some well known results, see \cite{suzC}. Some proofs are provided for reader's convenience.\\
The fact that for a given smooth initial value $\rho_0$, \eqref{system} is locally well posed in time is well known (see for example \cite{Bil}, \cite{GZ}). 
Also, the fact that $ \rho(\cdot,t) \geq 0$ follows for non negative initial value $\rho_0$ by the maximum principle. Actually it holds $ \rho(\cdot,t)  > 0$ for any 
$t>0$, and in particular solutions are classical, see for example Theorem 3.1 in \cite{suzC}. The total mass conservation for \eqref{system} can be obtained as follows:

\begin{equation}\label{mass}
\frac{d}{d t} \int_{\Omega} \rho \, dx = \int_\Omega \nabla \cdot (\nabla \rho - \rho \nabla (u + \log V)) \, dx = \int_{\partial \Omega} \left( \frac{\partial \rho}{\partial \nu} - \rho \frac{\partial (u + \log V)}{\partial \nu} \right) \, d\sigma = 0.
\end{equation}

Let us then consider the free energy associated to problem \eqref{system}, that is 

\begin{equation}\label{free}
\mathcal{F} (\rho) = \int_\Omega \rho  \left( \log \left( \frac{\rho}{V} \right) - 1 \right) - \frac{1}{2} \int_\Omega \rho \, 
G[\rho] = \int_\Omega \rho  \left( \log \left( \frac{\rho}{V} \right) - 1 - \frac{u}{2} \right), 
\end{equation}

where $G$ is the Green function of $-\Delta$ in $\Omega$ with homogeneous Dirichlet boundary conditions, $G[\rho](x)=(G \ast \rho) (x)$ 
and $\rho\in \mathcal{P}_\lambda$ where,
\begin{equation}\label{density}
\mathcal{P}_\lambda = \left\{ \rho \in L^1 (\Omega) : \rho \geq 0 \, \text{ a.e. }, \int_{\Omega} \rho = \lambda \, \text{ and } \, 
\int_{\Omega} \rho \log \rho < \infty \right\}.  
\end{equation}

Clearly, for any $t>0$ for which $\rho(x,t)$ is defined, we have $\rho\in \mathcal{P}_\lambda$ whenever $\rho_0$ is (say) smooth, see for example \cite{suzC}. 
Thanks to the conservation of mass \eqref{mass} we can deduce an important property of $P(\lambda,\Omega)$: the decrease of the free energy $\mathcal{F}$ along the flow associated with the first equation in \eqref{system}. In fact, notice that by \eqref{system} and integration by parts we have 

$$\int_\Omega \rho u_t = - \int_\Omega u_t \Delta u = \int_\Omega \nabla u_t \nabla u = \int_\Omega u \rho_t.$$

From this and the fact that $V$ is time independent, the homogeneous Neumann conditions and again integrating by parts, we deduce that:

$$\frac{d}{d t} \mathcal{F} (\rho) = \frac{d}{d t} \int_\Omega \rho  \left( \log \left( \frac{\rho}{V} \right) - 1 - \frac{u}{2} \right) \, dx = $$
$$=  \int_\Omega \rho_t  \left( \log \left( \frac{\rho}{V} \right) - \frac{u}{2} \right) \, dx - \int_\Omega \rho_t \, dx + \int_\Omega \rho \frac{d}{d t} \left( \log \left( \frac{\rho}{V} \right) - \frac{u}{2} \right) \, dx = $$
$$= \int_\Omega \rho_t  \left( \log \left( \frac{\rho}{V} \right) - \frac{u}{2} \right) \, dx + \int_\Omega \left( \rho_t - \frac{1}{2} \rho u_t\right) \, dx = $$
$$= \int_\Omega \rho_t  \left( \log \left( \frac{\rho}{V} \right) - u \right) \, dx =$$
$$= \int_\Omega \nabla \cdot (\nabla \rho - \rho \nabla (u + \log V)) \left( \log \left( \frac{\rho}{V} \right) - u \right) \, dx = $$
$$= \int_{\partial \Omega} \left( \frac{\partial \rho}{\partial \nu} - \rho \frac{\partial (u + \log V)}{\partial \nu} \right) \left( \log \left( \frac{\rho}{V} \right) - u\right) \, d\sigma$$
$$- \int_\Omega (\nabla \rho - \rho \nabla (u + \log V)) \nabla \left( \log \left( \frac{\rho}{V} \right) - u\right) \, dx =$$
$$- \int_\Omega \rho \left| \nabla \left( \log \left( \frac{\rho}{V} \right) - u\right) \right|^2 \, dx \leq 0.$$

Hence we get that:
\begin{equation}\label{dissip}
\mathcal{F} (\rho(\cdot,t)) \leq \mathcal{F} (\rho (\cdot,0)) = \mathcal{F} (\rho_0(\cdot)) \quad \forall \, t > 0.
\end{equation}

\bigskip
\bigskip

For any $v \in H^{1}_{0} (\Omega)$ let us consider the variational functional 
\begin{equation}\label{energy}
J_\lambda (v) = \frac{1}{2} \int_\Omega |\nabla v|^2 - \lambda \log \left( \int_\Omega V e^{v}\right). 
\end{equation}
It is well known that critical points $u_\lambda$ of $J_\lambda$ are weak solutions of \eqref{mf}. 
As mentioned in the introduction we will need the following recently derived \cite{BCV} existence result of
strict local minimizers for $J_\lambda$ with large mass. 
\brm If $v\in H^{1}_{0} (\Omega) $ then we define
$$
\|v\|_{H^{1}_{0} (\Omega) }=\|v\|_{L^2 (\Omega) }+\|\nabla v\|_{L^2 (\Omega) }.
$$
\erm

\bigskip

\begin{theorem}\label{subsup}

{\rm (a)} Let $\om$ be any smooth and simply connected domain.
For any $c\in(0,1]$ let $D = \frac{a}{b}$ (as given in \eqref{hypV}) and $c_{\sscp D} = c D$. Then there exist $\overline{\a}_*>\underline{\a}_*(c_{\sscp D})>0$ 
such that
if $\{\a^2x^2+y^2\leq \beta_-^2\}\subset\Omega\subset\{\a^2x^2+y^2\leq\beta_+^2\}$ with 
$c=\tfrac{\beta_-^2}{\beta_+^2}$ then, for any $\a\in(0,\underline{\a}_*(c_{\sscp D})]$ and for any $\lm\leq\lm_{\a,c_{\sscp D}}$, 
problem \rife{mf} admits a solution $u_{\lm}$. Here $\underline{\lm}_{\a,c_{\sscp D}}< \lm_{\a,c_{\sscp D}}<\overline{\lm}_\a$ and 
$\underline{\lm}_{\a,c_{\sscp D}}$, $\overline{\lm}_\a$ are strictly decreasing (as functions of $\a$) in $(0,\underline{\a}_*(c_{\sscp D})]$, $(0,\overline{\a}_*]$ 
respectively with $\underline{\lm}_{\underline{\a}_*(c_{\sscp D}),c_{\sscp D}}=8\pi=\overline{\lm}_{\overline{\a}_*}$ and 
$\underline{\lm}_{\a,c_{\sscp D}}\simeq\frac{4\pi c_{\sscp D}}{(8-c_{\sscp D})\a}$, $\ov{\lm}_{\a}\simeq\frac{2\pi}{3\a}$ as $\a\to 0^+$.\\

{\rm (b)}
There exists $\bar N> 4\pi$ such that if $\Omega$ is an open, bounded and convex set (therefore simple)
whose isoperimetric ratio, $N\equiv N(\Omega)=\frac{L^2(\pa\Omega)}{A(\Omega)}$, satisfies $N\geq\bar N$, then
for any $\lm\leq\lm_\textnormal{\tiny{$N$}}$ problem \rife{mf} admits a solution $u_{\lm}$. Here
$\underline{\Lambda}_\textnormal{\tiny{$N$},D}<\lambda_\textnormal{\tiny{$N$}}<
\overline{\Lambda}_\textnormal{\tiny{$N$}}$ with $\underline{\Lambda}_\textnormal{\tiny{$\bar N$},D}=
8\pi$, $\underline{\Lambda}_\textnormal{\tiny{$N$},D}$ and $\overline{\Lambda}_\textnormal{\tiny{$N$}}$
strictly increasing in $N$ and $\underline{\Lambda}_\textnormal{\tiny{$N$},D}\simeq\frac{D \pi^2 N}{16(32-D)}+\mbox{\rm O}(1)$,
$\overline{\Lambda}_\textnormal{\tiny{$N$}}\simeq\frac{2\sqrt3 N}{\pi}+\mbox{\rm O}(1)$ as $N\to+\infty$.\\

{\rm (c)} The solution $u_\lm$ found in both cases {\rm (a)} and {\rm (b)} is a strict local minimizer of $J_\lambda$ and 
the first eigenvalue of the linearized problem for \rife{mf} at $u_\lm$ is strictly positive. In particular, there exists $\delta_0>0$ such that 
\beq\label{220415.01}
J_\lambda(u)>J_\lambda(u_\lm)
\eeq
for any $0<\|u-u_\lm\|_{H^1_0(\om)}<\dt_0$. Moreover, there exists $\eta_0>0$ such that, for any $0<\eta_1<\eta_0$, there exists $d_1 > 0$ such that
\beq\label{220415.1}
J_\lambda(u)\geq J_\lambda(u_\lm) +d_1
\eeq
whenever $\eta_1\leq \|u-u_\lm\|_{L^2(\om)}<\eta_0$. 

\end{theorem}
\proof
We shall derive part {\rm (a)} of the statement and skip the details of the proof of part $(b)$ which can be handled by the same 
argument adopted in \cite{BdM2}.\\
In view of the dilation invariance of \eqref{mf} (that is, if $u_\lm(z)$ solves \rife{mf} in $\om$, then for any $t>0$ $u_\lm(tz)$ solves 
\rife{mf} in $(t)^{-1}\om$ with $V(z)$ replaced by $V(tz)$ which still satisfies \rife{hypV}) for fixed $c \in (0,1]$ and up to a rescaling we can assume without loss of generality that
$$\Omega_{\a,c}:=\{\a^2 x^2+y^2\leq c\}\subseteq\Omega\subseteq\{\a^2 x^2+y^2\leq 1\}=:\Omega_\a.$$
  
Let us consider the following Liouville-type \cite{Lio} problem
\begin{equation}\label{liou}
\begin{cases}
-\Delta u = \mu V e^u \quad\mbox{in }\Omega\\
\quad u = 0 \quad \mbox{on } \partial \Omega
\end{cases}
\end{equation}
and let us define
\beq\label{2}
v_{\a,\ga} = 2 \log{\left(\frac{(1+\ga^2)}{1+\ga^2(\a^2x^2+y^2)}\right)},\quad (x,y)\in \om_\a.
\eeq
A straightforward evaluation shows that $v_{\a,\ga}$ satisfies
\beq\label{3}
\graf{
-\Delta v_{\a,\ga}= V_{\a,\ga}{\dsp e^{v_{\a,\ga}}} & \mbox{in}\quad \om_\a\\
v_{\a,\ga}=0 & \mbox{on}\quad \pa\om_\a,
}
\eeq
where
\beq\label{4}
V_{\a,\ga}(x,y)= \frac{4 \ga^2}{(1+\ga^2)^2}\left(1+\a^2+\ga^2(1-\a^2)(\a^2x^2-y^2)\right).
\eeq
Since
$$V_{\a,\ga}(x,y)\geq g_{+}(\ga,\a):=\frac{4 \ga^2}{(1+\ga^2)^2}\left(1+\a^2+\ga^2(\a^2-1)\right),\quad \forall (x,y)\in \om_\a,$$
thanks to \eqref{hypV} we easily verify that $v_{\a,\ga}$ is a supersolution of \eqref{liou} whenever
\beq\label{5}
b \mu \leq g_{+}(\ga,\a).
\eeq

For fixed $\a \in (0,1)$, the function $h_\a(t)=g_{+}(\sqrt{t},\a)$ satisfies $h_\a(0)=0=h_\a\left(\frac{1+\a^2}{1-\a^2}\right)$, 
is strictly increasing in $\left(0,\frac{1+\a^2}{3-\a^2}\right)$ and strictly decreasing in $\left(\frac{1+\a^2}{3-\a^2},\frac{1+\a^2}{1-\a^2}\right)$.

Therefore, putting $\overline{\ga}_\a^2=\frac{1+\a^2}{3-\a^2}$ and 
$\overline{\mu}_{\a,b}:= \frac{1}{b} h_\a\left(\overline{\ga}_\a^2\right) \equiv \frac{1}{b} g_{+}(\overline{\ga}_\a,\a) \equiv \frac{(\a^2+1)^2}{2 b}$, 
we see in particular that for each $\mu\in(0,\overline{\mu}_{\a,b}]$ there exists a unique $\ga^{+}_\a\in\left(0,\overline{\ga}_\a\right]$ such that $\frac{1}{b} g_{+}(\ga^{+}_\a,\a)=\mu$ and $v_{\a,\ga^{+}_\a}$ is a supersolution of \eqref{liou}. Indeed we have

$$
\left(\ga^{+}_\a\right)^2=\left(\ga^{+}_\a(\mu)\right)^2=\frac{2(1+\a^2)- \mu b -2\sqrt{(1+\a^2)^2-2 \mu b }}{\mu b + 4(1-\a^2)}.
$$

On the other hand let us consider
\beq\label{vmeno}
v_{\a,\ga,c}=\left\{
               \begin{array}{ll}
                 2 \log{\left(\frac{(1+\ga^2)}{1+\tfrac{\ga^2}{c}(\a^2x^2+y^2)}\right)},& \hbox{$(x,y)\in\Omega_{\a,c}$} \\
                 0, & \hbox{$(x,y)\in\Omega\setminus\Omega_{\a,c}$.}
               \end{array}
             \right.
\eeq

Again a straightforward computation shows that $v_{\a,\ga,c}$ satisfies

$$
\left\{
  \begin{array}{ll}
-\Delta v_{\a,\ga,c}=V_{\a,\ga,c}{\dsp e^{v_{\a,\ga,c}}} & \mbox{in}\quad \Omega_{\a,c}\\
v_{\a,\ga,c}=0 & \mbox{on}\quad \pa\Omega_{\a,c},  \end{array}
\right.
$$

where

$$
V_{\a,\ga,c}(x,y)=\left\{
               \begin{array}{ll}
            \frac{4 \ga^2}{c(1+\ga^2)^2}\left(1+\a^2+\frac{\ga^2}{c}(1-\a^2)(\a^2x^2-y^2)\right) & \mbox{in }\Omega_{\a,c} \\
                 0 & \mbox{in }\Omega\setminus\Omega_{\a,c}.
               \end{array}
             \right.
$$

Since
$$
V_{\a,\ga,c}(x,y)\leq g_{-}(\ga,\a,c):=\frac{4 \ga^2}{c(1+\ga^2)^2}\left(1+\a^2+\ga^2(1-\a^2)\right),\quad \forall (x,y)\in \om,
$$

it is not difficult to check that $v_{\a,\ga,c}$, again by \eqref{hypV}, is a subsolution of \eqref{liou} whenever
\beq\label{6}
a \mu \geq g_{-}(\ga,\a,c).
\eeq

For fixed $\a\in(0,1)$, the function $f_{\a,c}(t)=g_{-}(\sqrt{t},\a,c)$, $t\in(0,\overline{\ga}_\a^2]$ is strictly increasing and
satisfies $f_{\a,c}(t)>h_\a(t)$. Once again, putting $\overline{\ga}_\a^2=\frac{1+\a^2}{3-\a^2}$ and 
$\underline{\mu}_{\a,a}:= \frac{1}{a} h_\a\left(\overline{\ga}_\a^2\right) \equiv \frac{1}{a} g_{-}(\overline{\ga}_\a,\a) \equiv \frac{(\a^2+1)^2}{2a}$, 
we see that for each $\mu\in(0,\underline{\mu}_{\a,a}]$ there exists a unique $\ga^{-}_{\a,c}\in\left(0,\overline{\ga}_\a\right)$ 
such that $\frac{1}{a} g_{-}(\ga^{-}_{\a,c},\a,c)=\mu$, $\ga^{-}_{\a,c}<\ga^{+}_\a$ and $v_{\a,\ga^{-}_{\a,c},c}$ is a subsolution of \eqref{liou}. 
Indeed we have

$$\left(\ga^{-}_{\a,c}\right)^2=\left(\ga^{-}_{\a,c}(\mu)\right)^2=\frac{\mu a c-2(1+\a^2)+2\sqrt{(1+\a^2)^2-2\a^2 \mu a c}}{4(1-\a^2)- \mu a c}.$$

Notice that since $a < b$ we clearly have $\underline{\mu}_{\a,a} > \overline{\mu}_{\a,b}$. In conclusion, since 
$\ga^{-}_{\a,c}(\mu)\leq \ga^{+}_\a(\mu)$ implies $v_{\a,\ga^{-}_{\a,c},c}\leq v_{\a, \ga^{+}_\a}$, 
for fixed $\a\in (0,1)$ and for each $\mu\in(0,\overline{\mu}_{\a,b}]$ we can set
$$\underline{u}_{\mu}=v_{\a,\ga^{-}_{\a,c}(\mu),c},\quad \overline{u}_\mu=v_{\a,\ga^{+}_\a(\mu)},$$
and conclude via well known sub-supersolution results \cite{clsw} that a solution (in a suitable weak sense) $u_{\a,\mu,c}$ of \eqref{liou} exists which satisfies
\beq\label{susu}
v_{\a,\ga^{-}_{\a,c}(\mu),c}\leq u_{\a,\mu,c}\leq v_{\a,\ga^{+}_\a(\mu)},\quad \forall (x,y)\in\om.
\eeq

Then the Brezis-Merle results \cite{bm} and standard elliptic regularity follows that $u_{\a,\mu,c}$ is a classical solution of \eqref{liou}.

Any such a solution $u_{\a,\mu,c}$ therefore solves \eqref{mf} with $\lm=\lm_{\a,c_{\sscp D}}(\mu)$ satisfying

\beq\label{stimalambdamu1}
\lm=\lm_{\a,c_D}(\mu)=\mu\int\limits_{\om} V e^{u_{\a,\mu,c}}\geq \mu a \int\limits_{\om_{\a,c}}e^{v_{\a,\ga^{-}_{\a,c}(\mu),c}}=
\mu a c\frac{\pi}{\a}\,(1+(\ga^{-}_{\a,c}(\mu))^2),
\eeq
and
\beq\label{stimalambdamu2}
\lm=\lm_{\a,c_D}(\mu)=\mu\int\limits_{\om} V e^{u_{\a,\mu,c}} \leq \mu b \int\limits_{\om_\a} e^{v_{\a,\ga^{+}_\a(\mu)}}=
\mu b \frac{\pi}{\a}\,(1+(\ga^{+}_\a(\mu))^2).
\eeq

In case $\mu=\overline{\mu}_{\a,b}$, recalling that $D = \frac{a}{b}$ and $c_{\sscp D} = c D$, we have 
$$(\ga^{-}_{\a,c}(\overline{\mu}_\a))^2\equiv
\underline{\ga}_{\a,c}^2 = (1+\a^2) \frac{c_{\sscp D} (1+\a^2) + 4 \left( \sqrt{1- c_{\sscp D} \a^2} -1 \right)}{8(1-\a^2)- c_{\sscp D} (1+\a^2)^2}$$
and
$$(\ga^{+}_\a(\overline{\mu}_\a))^2 \equiv \overline{\ga}_\a^2 = (1+\a^2) \frac{3-\a^2}{8(1-\a^2)+(1+\a^2)^2},$$

so that, by \eqref{stimalambdamu1} and \eqref{stimalambdamu2} respectively we have:

\bel\label{lambdasotto}
\lm_{\a,c_{\sscp D}}:=\lm_{\a,c_{\sscp D}}(\overline{\mu}_{\a,b}) \geq \underline{\lm}_{\a,c_{\sscp D}}=\frac{a c(1+\a^2)^2}{2 b} \frac{\pi}{\a}(1+\underline{\ga}_{\a,c}^2)\simeq
\frac{4\pi c_{\sscp D}}{(8-c_{\sscp D})\a},
\eel
and
\bel\label{lambdasopra}
\lm_{\a,c_{\sscp D}}\leq \overline{\lm}_{\a}=\frac{(1+\a^2)^2}{2}\frac{\pi}{\a}
(1+\overline{\ga}_\a^2)\simeq \frac{2 \pi}{3\a}
\eel
as $\a\to 0^+$. Moreover it is easy to verify that $\underline{\lm}_{\a,c_{\sscp D}}$ is strictly decreasing at least
for $\a\in(0,\frac 1{2\sqrt{10}}]$ and that there exists
$\underline{\a}_*(c_{\sscp D})<\frac 1{2\sqrt{10}}$ such that
$\underline{\lm}_{\a,c_{\sscp D}}\geq 8\pi$ for any $\a\in(0,\underline{\a}_*(c_{\sscp D})]$. These estimates are uniform in
 $0<D\leq 1$ and in $c\in(0,1]$. We also see that $\overline{\lm}_\a \to \left(\frac{4\pi}{b}\right)^{-}$ as $\a\to 1^{-}$,
is strictly decreasing for $\a\in(0,\a_p]$ and strictly increasing
for $\a\in[\a_p,1)$ for some $\a_p\simeq 0.5$ and then it is straightforward to check that
there exists $\overline{\a}_*>\underline{\a}_*(c)$  such that
$\overline{\lm}_\a\geq 8\pi$ for any $\a\in(0,\overline{\a}_*]$. Finally, since $\lm_{\a,c_D}(\mu)$ is continuous in $\mu$ and by using \eqref{stimalambdamu1} and \eqref{stimalambdamu2}
$$
0<\lm_{\a,c_D}(\mu)\leq\mu\frac{\pi}{\a}\,(1+(\ga^{+}_\a(\mu))^2)\stackrel{\textnormal{as $\mu\to0$}}{\longrightarrow} 0,
$$
we obtain the existence of a solution for $P(\lambda,\Omega)$ not only for $\lambda=\lm_{\a,c_D}$, but for any $\lambda\in(0,\lm_{\a,c_D}]$ as well. This fact concludes the proof of the existence result claimed in {\rm (a)}.\\

Next we prove that the solutions obtained so far are strict local minimizers of $J_{\lm}$, i.e. part $(c)$ of the statement. Actually we need a stronger result, that is, the linearized 
problem relative to \eqref{mf} has a strictly positive first eigenvalue. 
Putting
$$
\de=\de(u)=\dfrac{V e^u}{\int\limits_{\om} V e^u},\quad \mbox{and}\quad <f>_\de=\int\limits_{\om}\de(u)f,
$$
then the linearized problem for \eqref{mf} takes the form
\bel\label{7.1}
\left\{
  \begin{array}{ll}
     -\Delta \vp - \lm\de(u)\vp +\lm \de(u)<\vp>_\de =0 & \mbox{in}\ \ \om \\
    \vp =0 & \mbox{on}\ \ \pa\om.
  \end{array}
\right.
\eel
Letting $H\equiv H^{1}_0(\om)$ and
$$
\mathcal{L}(\phi,\psi)= \int\limits_{\om} \left(\nabla \phi\cdot\nabla \psi\right)-
\lm\int\limits_{\om} \de(u)\phi\psi +
\lm\left(\int\limits_{\om} \de(u)\phi\right)\left(\int\limits_{\om} \de(u)\psi\right),\;(\phi,\psi)\in H\times H,
$$
then by definition $\vp\in H$ is a weak solution of \eqref{7.1} if
$$
\mathcal{L}(\vp,\psi)=0,\quad \forall\, \psi \in H.
$$
We define $\tau\in\erre$ to be an eigenvalue of the operator
$$
L[\vp]:=-\Delta\vp  - \lm\de(u) (\vp -< \vp >_\de),\quad \vp\in H,
$$
if there exists a weak solution $\phi_0\in H\setminus \{0\}$ of the linear problem
\beq\label{7.2}
 -\Delta \phi_0 - \lm\de(u)\phi_0+\lm \de(u)<\phi_0>_\de =\tau \de(u)\phi_0\ \ \mbox{in}\ \ \om,
\eeq
that is, if
$$
\mathcal{L}(\phi_0,\psi)=\tau\int\limits_{\om} \de(u)\phi_0\psi,\quad \forall\, \psi \in H.
$$

Standard arguments show that the eigenvalues form an unbounded (from above) sequence
$$
\tau_1\leq \tau_2\leq\cdots\leq \tau_n\cdots,
$$

with finite dimensional eigenspaces  (although the first eigenfunction changes sign and cannot be assumed to be 
simple in this situation).

Let us define
$$
Q(\phi)=\frac{\mathcal{L}(\phi,\phi)}{<\phi^2>_\de}=
\frac{\int\limits_{\om} \left|\nabla \phi\right|^2-\lm<\phi^2>_\de+\lm <\phi>_\de^2}{<\phi^2>_\de},\quad \phi\in H.
$$
 In particular it is not difficult to prove that the first eigenvalue can be characterized as
follows
$$
\tau_1=\inf\{Q(\phi)\,|\, \phi\in H\setminus\{0\}\}.
$$

At this point we argue by contradiction and assume that $\tau_1\leq 0$. Therefore we readily conclude that
$$
\tau_0:=\inf\{Q_0(\phi)\,|\, \phi\in H\setminus\{0\}\}\leq 0,\quad\mbox{where}\quad
Q_0(\phi)=\frac{\mathcal{L}_0(\phi,\phi)}{<\phi^2>_\de}
$$
and
$$
\mathcal{L}_0(\phi,\psi)= \int\limits_{\om} \left(\nabla \phi\cdot\nabla \psi\right)
-\lm\int\limits_{\om} \de(u)\phi\psi,\;(\phi,\psi)\in H\times H.
$$

Clearly $\tau_0$ is attained by a simple and positive eigenfunction $\vp_0$ which satisfies
\bel\label{7}
\left\{
  \begin{array}{ll}
    \dsp -\Delta\vp_0 -\lm \de(u)\vp_0 =\tau_0 \de(u)\vp_0 & \mbox{in}\ \ {\om} \\
    \vp_0 =0 & \mbox{on}\ \ \pa{\om}.
  \end{array}
\right.
\eel

Let us recall that we have obtained solutions for \eqref{mf} as solutions of \rife{liou}
in the form $u=u_{\a,{\mu},c}$, for some $\mu=\mu(\a)\leq \ov{\mu}_{\a,b}$ whose value of $\lm=\lm(\mu,\a,c_{\sscp D})$
was then estimated as a
function of $\a$. That point of view is well suited for our purpose,  that is, we get back to $\mu=\lm\left(\int_{\om} V e^u\right)^{-1}$.
Hence, let us observe that for a generic value $\mu\leq \ov{\mu}_{\a,b}$ \rife{7} takes the form
\bel\label{9}
\left\{
  \begin{array}{ll}
    \dsp -\Delta\vp_0 -\mu V K_{\a,\mu,c} \vp_0  =\nu_0 V K_{\a,\mu,c} \vp_0 & \mbox{in}\ \ {\om} \\
    \vp_0 =0 & \mbox{on}\ \ \pa{\om},
  \end{array}
\right.
\eel
where
$$
K_{\a,\mu,c}=e^{u_{\a,\mu,c}}\quad\mbox{and}\quad \nu_0=\mu\frac{\tau_0}{\lm}\leq 0.
$$

We observe that, by defining
$$
K_{\a,\mu,c}^{(-)}:=e^{v_{\a,\ga^{-}_{\a,c}(\mu),c}}=\left\{
                                                       \begin{array}{ll}
                                                         \left(\frac{1+\ga^{-}_{\a,c}(\mu)^2}{1+\frac{\ga^{-}_{\a,c}(\mu)^2}{c}(\a^2x^2+y^2)}\right)^2 & \mbox{$(x,y)\in\Omega_{\a,c}$} \\
                                                         1 & \mbox{$(x,y)\in\Omega\setminus\Omega_{\a,c}$,}
                                                       \end{array}
                                                     \right.
$$
$$
K_{\a,\mu}^{(+)}:=e^{v_{\a,\ga^{+}_\a(\mu)}}=\left(\frac{1+\ga^{+}_\a(\mu)^2}{1+\ga^{+}_\a(\mu)^2(\a^2x^2+y^2)}\right)^2,\qquad (x,y)\in\Omega_{\a}
$$
we have
$$
K_{\a,\mu,c}^{(-)}\leq K_{\a,\mu,c}\leq K_{\a,\mu}^{(+)}\qquad \textnormal{for any }(x,y)\in\Omega.
$$

In particular, since
$$
K_{\a,\mu}^{(+)}\leq (1+\ga_\a^{+}(\mu)^2)^2\quad\mbox{and}\quad
1\leq K_{\a,\mu,c}^{(-)}\leq (1+\ga_{\a,c}^{-}(\mu)^2)^2\quad \mbox{in}\quad \om,
$$
and
\beq\label{140213.0}
\om\subset T_\a:=\{(x,y)\in\rdue\,|\,|\,x|\leq (\a)^{-1},\;|\,y|\leq 1\},
\eeq

then, by using the fact that
$$
\nu_0=\inf\left\{
\left.\frac{\int\limits_{\om} \left|\nabla \vp\right|^2\,dx-\mu\int\limits_{\om} K_{\a,\mu} \vp^2\,dx}
{\int\limits_{\om} K_{\a,\mu}\vp^2\,dx}\;\right|
\,\vp\in H\right\}\leq0,
$$
it is not difficult to check that, for some $\mu\leq \ov{\mu}_{\a,b}=\frac{(1+\a^2)^2}{2b}$, thanks to \eqref{hypV} the following inequality
holds:
\beq\label{13}
\inf\left\{
\left.\frac{\int\limits_{T_\a} \left|\nabla \vp\right|^2\,dx-\mu b (1+\ga_\a^{+}(\mu)^2)^2\int\limits_{T_\a} \vp^2\,dx}
{\int\limits_{T_\a} \vp^2\,dx}\;\right| \,\vp\in H\right\}\leq 0.
\eeq

Hence, there exists $\overline{\mu}_0\leq 0$ such that, putting $\sg=\sg(\mu,\a)=\mu b (1+\ga_\a^{+}(\mu)^2)^2+\overline{\mu}_0$, there exists a weak solution $\phi_0\in H$ of
\bel\label{14}
\left\{
  \begin{array}{ll}
    \dsp -\Delta\phi_0 -\sg\phi_0 =0 & \mbox{in}\ \ T_\a,\\
    \dsp \phi_0=0 & \mbox{on}\quad  \pa T_\a.
  \end{array}
\right.
\eel
It is well known that the minimal eigenvalue $\sg_{\it min}$ of \eqref{14} satisfies
$\sg_{\it min}=\frac{\pi^2}{4}\a^2+\frac{\pi^2}{4}>2(1+\a^2)$ and we conclude that
\beq\label{16}
2(1+\a^2)\leq\sg(\mu,\a)=\mu b (1+\ga_\a^{+}(\mu)^2)^2+\overline{\mu}_0.
\eeq
Next, since $\underline{\a}_*(c_{\sscp D})<\tfrac{1}{2\sqrt{10}}$, it is not difficult to check that
$\sg=\sg(\mu,\a)$ satisfies
$$
\sg(\mu,\a)\leq 1,
$$
for any $\a\leq\underline{\a}_*(c_{\sscp D})$, which is of course a contradiction to \eqref{16}. Therefore $\tau_1$ is strictly positive as claimed, and this clearly yields \eqref{220415.01} for some small $\delta_0 > 0$.\\

Finally we prove \rife{220415.1}: on the one hand, since $\|u-u_\lm\|_{H^1_0(\om)} < \delta_0$, then by Poincar\'e inequality with constant $\mu_1 = \mu_1 (\Omega)$ we have $\|u-u_\lm\|_{L^2 (\om)} < \eta_0 := \frac{\delta_0}{\mu_1}$. On the other hand, by a simple Taylor expansion for $J_\lm$ around the local minimizer $u_\lm$ and by choosing $\eta_0$ smaller if necessary, we see that
$$J_\lambda(u) = J_\lambda(u_\lm) + \frac{1}{2} <J_\lambda ''(u_\lambda) [u-u_\lambda],u-u_\lambda> + \mbox{\rm o}(\|u-u_\lm\|_{L^2 (\om)}^{2})$$ 
$$\geq J_\lambda(u_\lm) + \frac{1}{4} <J_\lambda ''(u_\lambda) [u-u_\lambda],u-u_\lambda>$$
for $\|u-u_\lm\|_{L^2(\om)}<\eta_0$. Now, thanks to the fact that the linearized operator $\mathcal{L}$ admits a positive first eigenvalue $\tau_1$ we have:
$$<J_\lambda ''(u_\lambda) [u-u_\lambda],u-u_\lambda> = \mathcal{L} (u-u_\lm, u-u_\lm) \geq \tau_1 <(u-u_\lm)^2>_\de$$
Notice that, thanks to \eqref{hypV} and \eqref{susu}, there exists $c_\lm>0$ such that $\omega(u_\lm) \geq c_\lm$, which spells that
$$<(u-u_\lm)^2>_\de \geq c_\lm \|u-u_\lm\|_{L^2(\om)}$$
Thus, we finally find that for any $\eta_1 < \eta_0$ we have
$$J_\lambda(u) \geq J_\lambda(u_\lm) + \frac{\tau_1 c_\lm \eta_1}{4}$$
whenever $\eta_1 \leq \|u-u_\lm\|_{L^2 (\om)}< \eta_0$, which proves \eqref{220415.1} with $d_1 := \frac{\tau_1 c_\lm \eta_1}{4}$.
\fineproof

\bigskip
\bigskip

\section{The existence of free energy minimizers with large masses.}\label{s3}

As mentioned in the introduction we will need some properties in the dual Orlicz topology that the free energy minimizers inherit from those of the $u_\lm$ found in Theorem \ref{subsup}. In order to clarify this aspect let 
$$
L_{\Phi} (\om) := \{\,\rho\,\mbox{measurable in}\,\,\om \,:\,\| \rho \|_{\Phi} < +\infty \},
$$

be the Orlicz space \cite{Ad}, \cite{KJF} whose Young functions are
$$
\Phi(t)=t\log{(1+t)},\quad t\geq 0,\quad \Psi(s)=\max\limits_{t\geq 0}\{ts-\Phi(t)\},\quad s\geq 0,
$$
where
$$
\|\rho\|_{\Phi}:= \sup\limits_{h}\left\{\left|\;\ino \rho h\;\;\right|,\;\ino \Psi(|h|)\leq 1 \right\}.
$$

The Orlicz space contains the so called Orlicz class of all functions $\rho$ which are measurable in $\om$ and satisfy $\ino \Phi(|\rho|)<+\infty$. 
For any $u$ which satisfies $\ino \Psi(|u|)<+\infty$, let us also introduce the Luxemburg norm, 
\beq\label{dualnorm}
\|u\|_{\Psi}:=\inf\left\{\alpha >0\,:\, \ino \Psi(\alpha^{-1}|u|)\leq 1\right\}.
\eeq
It is straightforward to check that
\beq\label{130613.1}
\Psi(s)\leq se^{s-1},\,s\in [0,+\infty)\quad\mbox{and} \quad \Psi(s)\leq e^{s-1}-(s-1),\;s\geq 2,
\eeq
whence in particular the Moser-Trudinger inequality \cite{moser} shows that \rife{dualnorm} is well defined for any $u\in H^{1}_0(\om)$.
It is well known \cite{Ad}, \cite{KJF} that $L_{\Phi}(\om)$ is a Banach space with respect to the $\|\cdot\|_{\Phi}$-norm. In particular the following version of the H\"{o}lder inequality holds
\beq\label{HoldOrl}
\ino\rho u\leq \|\rho\|_\Phi \|u\|_\Psi,
\eeq
for any $\rho\in L_{\Phi}(\om)$ and any $u$ such that \rife{dualnorm} is well defined.
In particular we  have
\beq\label{l1l2}
 \frac12 t^2\leq \Phi(t)\leq t^2,\,t\leq 1,\quad \mbox{and}\quad t\leq \Phi(t)\leq t^2,\,t\geq 1.
\eeq
Clearly $\mathcal{P}_{\lambda}$ is a convex subset of $L_{\Phi}(\om)$.  We will need the following result about 
the continuity of $\mathcal{F}$ with respect to the topology induced by $L_{\Phi}(\om)$. Let 
$$
\mathcal{P}:=
\left\{ \rho \in L^1 (\Omega) : \rho \geq 0 \, \text{ a.e. } \, \text{ and } \, 
\int_{\Omega} \rho \log \rho < \infty \right\}.
$$
For any density $\rho \in \mathcal{P}$ we will let $u_{\rho}$ 
be the corresponding potential, that is $u_{\rho}(x) = G[\rho](x)=(G \ast \rho) (x)$. Clearly  $u_{\rho}$ is the unique distributional solution of:
\begin{equation}\label{conc}
\begin{cases}
-\Delta u_{\rho} = \rho \quad\mbox{in }\Omega\\
\quad u_{\rho} = 0 \quad \mbox{on } \partial \Omega
\end{cases}
\end{equation}

\brm\label{remh1}
By using the Green's representation formula
\beq\label{Green}
u_{\rho} (z)=\ino G(z,w) \rho (w) \,dw,\quad \fo z \in \om,
\eeq
and the H\"{o}lder's inequality \rife{HoldOrl} we see that $u_\rho\in L^{\infty}(\om)$. 
Indeed, for any $\al \geq \frac{1}{2\pi}$ we have
$$
\ino \Psi\left(\frac{|G(z,w)|}{\al}\right)= \ino \Psi\left(\frac{G(z,w)}{\al}\right)dw \leq 
\ino \left(\textstyle \frac{G(z,w)}{\al}\right) e^{\sscp\frac{G(z,w)}{\al}}dw \leq
$$
$$
\frac{C_{\om}}{\al} \int\limits_{B_1(z)}\log\frac{1}{|z-w|}\left(\frac{1}{|z-w|} \right)^{\sscp \frac{1}{2\pi\al}}dw \leq 
\frac{C_{\om}}{\al} \int\limits_{B_1(z)}\log\frac{1}{|z-w|}\left(\frac{1}{|z-w|} \right)dw\leq \frac{C_0}{\al},
$$
so that 
\beq\label{200415.3}
\sup\limits_{z\in\om}\|G(z,\cdot)\|_{\Psi}\leq C_0
\eeq
for some constant $C_0\geq 1$ depending only by $\om$. Therefore, in particular $\rho u_\rho\in L^1(\om)$ and then standard 
truncation arguments show that $u_\rho\in H^1_0(\om)$. We conclude that $u_\rho$ is also a weak solution of \rife{conc}.
\erm

Then we have
\bpr
The functional $\mathcal{F}$ is continuous on $\mathcal{P}$ 
with respect to the topology induced by $L_{\Phi}(\om)$.
\epr
\proof
Fix $\rho \in \mathcal{P}$ and let $\{\rho_n\}\subset  \mathcal{P}$ be any sequence such that
\beq\label{190415.1}
\|\rho_n-\rho\|_{\Phi}\to 0,\,\ainf.
\eeq
We are going to prove that $\mathcal{F}(\rho_n)\to \mathcal{F}(\rho)$. 
We recall that, since $\Phi$ satisfies the following $\varDelta_2$-condition
\beq\label{delta2}
\Phi(2t)\leq 4 \Phi(t),\;\fo\,t\geq 1
\eeq
then, see \cite{KJF}, a sequence $\{\rho_n\}\subset L_{\Phi}(\om)$ 
satisfies \rife{190415.1} if and only if  
\beq\label{170415.1}
\lim\limits_{n\to +\infty} \ino \Phi(|\rho_n-\rho|)=0.
\eeq
Clearly $\rho_n$ satisfies \rife{170415.1} since $\mathcal{P}\subset L_{\Phi}(\om)$. By using \rife{l1l2} we find
$$
\ino |\rho_n-\rho| \leq \int\limits_{|\rho_n-\rho|\leq 1}|\rho_n-\rho|+\int\limits_{|\rho_n-\rho|\geq 1}|\rho_n-\rho|\leq
$$
$$
|\om|^{\frac12}\left(\,\int\limits_{|\rho_n-\rho|\leq 1}|\rho_n-\rho|^2\right)^{\frac12}+\int\limits_{|\rho_n-\rho|\geq 1}|\rho_n-\rho|\leq
$$
$$
|\om|^{\frac12}\left(\,\ino 2\Phi(\rho_n-\rho)\right)^{\frac12}+\ino \Phi(\rho_n-\rho).
$$

Therefore $\rho_n\to \rho$ in $L^1(\om)$. Next observe that
$$
u_{n}(z):=G[\rho_n](z)=\ino G(z,w) \rho_n(w), \quad \mbox{and} \quad u_\rho(z)=G[\rho](z)=\ino G(z,w) \rho(w),
$$
satisfy $u_n\to u_\rho$ in $L^{\infty}$. In fact, letting $z_n$ be any point where the maximum of $|u_n-u_\rho|$ is attained, we find
$$
\|u_n-u_\rho\|_\i=|u_n(z_n)-u_\rho(z_n)|\leq \ino G(z_n,w)|\rho_n(w)-\rho(w)|\leq \|\rho_n-\rho\|_{\Phi} \|G(z_n,\cdot)\|_{\Psi}\leq C_0 \|\rho_n-\rho\|_{\Phi},
$$
where we used \rife{HoldOrl} and \rife{200415.3}. Therefore, since $\rho_n$ converges in $L^1(\om)$ to $\rho$, 
we conclude that
$$
\ino \rho_n G[\rho_n]=\ino \rho_n u_n \to \ino \rho u_\rho=\ino \rho G[\rho],\ainf,
$$
by the duality $L^1(\om), L^{\i}(\om)$. This fact shows that the functional $\rho\mapsto \ino \rho G[\rho]$ is continuous.\\

Since we have shown that $\ino \rho G[\rho]$ is continuous and that $\rho_n\to \rho$ in $L^1(\om)$, then, 
to conclude the proof, we just need to show that in fact $\ino \rho_n\log(\rho_n) \to \ino \rho \log(\rho)$.  To this aim we observe that since $\rho_n$ 
satisfies \rife{190415.1}, and since $L_\Phi(\om)$ is a Banach space, then $\|\rho_n\|_{\Phi}$ is uniformly bounded. But then, see \cite{KJF} \S 3.10.9, 
since $\Phi$ satisfies \rife{delta2}, then there exists $C_2>0$ depending only by $\ino \rho\log\rho$ such that 
$$
\ino \Phi (\rho_n)\leq C_2,\;\fo  n\in \mathbb{N},
$$
whence in particular 
\beq\label{190415.2}
\ino \rho_n\log(\rho_n)\leq 
\int\limits_{\rho_n\geq 1} \rho_n\log(\rho_n) \leq\int\limits_{\rho_n\geq 1} \rho_n\log(1+\rho_n) \leq C_2,\;\fo n\in\mathbb{N}.
\eeq
At this point, for any $\eps>0$ we can choose $m_\eps\in (0,1)$ such that, setting 
$$
\om_{m,1}:=\{ \rho_n\leq m \} \cup \{\rho \leq m\},
$$
then, for any $m<m_\eps$ it holds
\beq\label{c1}
\left|\, \int\limits_{\om_{m,1}}  \rho_n \log(\rho_n) - \int\limits_{\om_{m,1}}  \rho \log(\rho) \right|\leq 
 \int\limits_{\om}  2m |\log(m)|<\eps, 
\eeq
where we used the fact that $m\log(m)\to 0$ as $m\to 0^+$. Next let us set $\om_m:=\om\setminus \om_{m,1}$ and decompose
\beq\label{c1.1}
\left|\ino \rho_n\log(\rho_n)-\ino \rho\log(\rho)\right|\leq 
\eeq
$$
\left|\,\int\limits_{\om_{m,1}} \rho_n\log(\rho_n)- \int\limits_{\om_{m,1}}\rho\log(\rho)\,\right|+
\left|\,\int\limits_{\om_m} \rho_n\log(\rho_n)-\int\limits_{\om_m}\rho\log(\rho)\,\right|\leq
$$
$$
2\eps +\left|\,\int\limits_{\om_m\cap \{\rho_n\geq \rho \}}  \rho_n\log(\rho_n)-\int\limits_{\om_m\cap \{\rho_n\geq \rho\}} \rho\log(\rho)\,\right|+
\left|\,\int\limits_{\om_m\cap \{\rho_n\leq \rho\}} \rho_n\log(\rho_n)-\int\limits_{\om_m\cap \{\rho_n\leq \rho\}}\rho\log(\rho)\,\right|
$$

We will use fact that for any $\al\geq 1$, we have

$$
\ino \Psi\left(\frac{|\log(\rho_n)|\chi_{\sscp \om_m}}{\al}\right)\leq \ino 
\left(\textstyle \frac{|\log(\rho_n)|\chi_{\sscp \om_m}}{\al}\right) e^{\sscp\frac{|\log(\rho_n)|\chi_{\sscp \om_m}}{\al} }=
$$
$$
\int\limits_{\om_m}\left(\textstyle \frac{|\log(\rho_n)|}{\al}\right) e^{\sscp\frac{|\log(\rho_n)|}{\al} }\leq 
\int\limits_{m\leq \rho_n \leq 1}\left(\textstyle \frac{|\log(\rho_n)|}{\al}\right) e^{\sscp\frac{|\log(\rho_n)|}{\al} }+
\int\limits_{\rho_n\geq 1}\left(\textstyle \frac{|\log(\rho_n)|}{\al}\right) e^{\sscp\frac{|\log(\rho_n)|}{\al} }\leq
$$
$$
|\om| \frac{|\log(m)|}{\al}\left(\frac{1}{m}\right)^{\frac{1}{\al}}+\int\limits_{\rho_n\geq 1} \frac{\log(\rho_n)}{\al}(\rho_n)^{\frac{1}{\al}}\leq 
|\om| \frac{|\log(m)|}{\al}\left(\frac{1}{m}\right)+\int\limits_{\rho_n\geq 1} \frac{\log(\rho_n)}{\al}(\rho_n)=
$$
$$
\frac{1}{\al}\left(\frac{|\om||\log(m)|}{m}+ \int\limits_{\rho_n\geq 1} \rho_n\log(\rho_n)\right)\leq 
\frac{1}{\al}\left(\frac{|\om||\log(m)|}{m}+ C_2\right),
$$
showing that, in view of \rife{190415.2}, 
$$
\|\log(\rho_n)\chi_{\sscp \om_m}\|_{\Psi}\leq C_*,
$$
for some constant $C_*$ depending only by $m$, $|\om|$ and $C_2$. Therefore we can estimate 
\beq\label{c1.2}
\left|\,\int\limits_{\om_m\cap \{\rho_n\geq \rho \}} \rho_n\log(\rho_n)-\int\limits_{\om_m\cap \{\rho_n\geq \rho \}}\rho\log(\rho)\,\right|\leq
\eeq
$$
\int\limits_{\om_m\cap \{\rho_n\geq \rho \}} |\rho_n-\rho||\log(\rho_n)|+\int\limits_{\om_m\cap \{\rho_n\geq \rho \}}\rho|\log(\rho_n)-\log(\rho)|\leq 
$$
$$
\|\rho_n-\rho\|_{\Phi}\|\log(\rho_n)\chi_{\sscp \om_m}\|_{\Psi} +\int\limits_{\om_m\cap \{\rho_n\geq \rho \}}\rho\frac{1}{\min\{\rho_n,\rho\}}|\rho_n-\rho|\leq 
$$
$$
C_*\|\rho_n-\rho\|_{\Phi}+\int\limits_{\om_m\cap \{\rho_n\geq \rho \}}\rho\frac{1}{\rho}|\rho_n-\rho|\leq
$$
$$
C_*\|\rho_n-\rho\|_{\Phi}+\int\limits_{\om}|\rho_n-\rho|<\eps,
$$
for any $n$ large enough, possibly  depending on $m$ and $C_2$, where
we used the mean value theorem. The same argument in a slightly easier form shows that  
\beq\label{c1.3}
\left|\,\int\limits_{\om_m\cap \{ \rho_n\leq \rho \}} \rho_n\log(\rho_n)-\int\limits_{\om_m\cap \{\rho_n\leq \rho\}}\rho\log(\rho)\,\right|<\eps
\eeq
for any $n$ large enough, possibly  depending on $m$, and we skip the details relative to this estimate to avoid repetitions. 
The fact that $\ino \rho_n \log(\rho_n) \to \ino \rho \log(\rho)$ is an immediate consequence of \rife{c1.1}, \rife{c1.2} and \rife{c1.3}.\fineproof

\bigskip
\bigskip

Next we have
\bpr\label{min}
Let $u_\lambda$ be a strict local minimum of $J_\lambda$ and assume that \rife{220415.01}, \rife{220415.1} hold.
If 
\beq\label{030415.2}
\rho_{0,\lm}= \lambda \dfrac{ V e^{\dsp u_\lambda}}{\int_\Omega V e^{\dsp u_\lambda}},
\eeq
then $\rho_{0,\lm}\in\mathcal{P}_{\lm}$ and the following property ${\bf (H)_{\lm}}$ holds:
there exist  and $\varepsilon_0 >\varepsilon_1> 0$ such that 
\beq\label{200415.01}
\mathcal{F} (\rho) - \mathcal{F} (\rho_{0,\lm}) > 0,
\eeq
for any $\rho\in \mathcal{P}_{\lm}$ such that $0 < \| \rho - \rho_{0,\lm} \|_{\Phi} < \varepsilon_0$ 
and
\beq\label{200415.1}
\mathcal{F} (\rho) \geq \mathcal{F} (\rho_{0,\lm}) + d_1,
\eeq
for any $\rho\in L^2(\om)\cap \mathcal{P}_\lm$ such that $\| \rho - \rho_{0,\lm} \|_{\Phi} = \varepsilon_1$, where $d_1$ is the constant 
introduced in \rife{220415.1}.
\epr
\proof
Whenever $u_{\rho}= G[\rho](x)=(G \ast \rho) (x)$ we set 
\beq\label{sgrho}
\sg_{u_{\rho}} := \lambda \dfrac{ V e^{\dsp u_{\rho}}}{\int_\Omega V e^{\dsp u_{\rho}}}.
\eeq 

If  $\rho_{0,\lm}$ is defined as in \rife{030415.2}, then $\rho_{0,\lm}\in \mathcal{P}_\lm$ and for any $\rho\in \mathcal{P}_\lm$ we find,
$$\mathcal{F} (\rho) - \mathcal{F} (\rho_\lambda) =$$
$$= \int_\Omega \rho  \left( \log \left( \frac{\rho}{V} \right) - 1 - \frac{1}{2} G[\rho] \right) 
- \int_\Omega \rho_{0,\lm}  \left( \log \left( \frac{\rho_{0,\lm}}{V} \right) - 1 - \frac{1}{2} G[\rho_{0,\lm}] \right) =$$
$$= \int_\Omega \rho  \left( \log \left( \frac{\rho}{V} \right) - \frac{1}{2} G[\rho] \right) + 
\lambda \log \left( \int_\Omega V e^{u_\lambda}\right) - \frac{1}{2} \int_\Omega |\nabla u_\lambda |^2 - 
\lambda \log \lambda - \int_\Omega \rho + \int_\Omega \rho_{0,\lm}=$$
$$= \int_\Omega \rho \log \left( \frac{\rho}{\sg_{u_{\rho}} } \right) + \int_\Omega \rho \log \left( \frac{\sg_{u_{\rho}} }{V} \right) 
- \frac{1}{2} \int_\Omega \rho   G[\rho] - \lambda \log \lambda - J_\lambda (u_\lambda) \geq$$
$$= \int_\Omega \rho \log \left( \frac{\sg_{u_{\rho}} }{V} \right) - \frac{1}{2} \int_\Omega \rho  G[\rho] - \lambda \log \lambda - J_\lambda (u_\lambda) =$$
$$= \int_\Omega \rho u_{\rho} - \frac{1}{2} \int_\Omega \rho  G[\rho] - \lambda \log \left( \int_\Omega V e^{u_\rho}\right) - J_\lambda (u_\lambda) =$$
\beq\label{030415.1}= J_\lambda (u_{\rho}) - J_\lambda (u_\lambda),
\eeq

where we have used Remark \ref{remh1}, \rife{sgrho} and the following facts: 

\begin{itemize}
\item by definition $\log \left( \frac{\rho_{0,\lm}}{V} \right) = \log \left(  \frac{\lambda e^{u_\lambda}}{\int_\Omega V e^{u_\lambda}} \right) = 
u_\lambda + \log \lambda - \log \left( \int_\Omega V e^{u_\lambda}\right) $;

\item by \eqref{mf} and the definition of $\rho_{0,\lm}$ we have: $-\int_\Omega \rho_{0,\lm} u_\lambda + 
\frac{1}{2} \int_\Omega \rho_{0,\lm} G [\rho_{0,\lm}] = -\frac{1}{2} \int_\Omega |\nabla u_\lambda |^2$;

\item since $\{\sg_{u_{\rho}} , \rho_{0,\lm},\rho \} \subset \mathcal{P}_\lambda$ we have 
$\ino \sg_{u_{\rho}} = \ino \rho_{0,\lm} = \ino \rho=\lambda$;

\item by Jensen's inequality applied to $\phi(t) = t \log t$ and $t = \frac{\rho}{\sg_{u_{\rho}}}$ we have:
$\int_\Omega \rho \log \left( \frac{\rho}{\sg_{u_{\rho}}} \right) \geq 0$.
\end{itemize}
We learned of this nice application of the Jensen's inequality in \cite{w2}. In view of \rife{220415.01} we have 
$J_\lambda (u) - J_\lambda (u_\lambda) > 0$ for any
$0 < \| u - u_{\lambda} \|_{H^{1}_{0} (\Omega)} < \delta_0$.  Therefore, to prove \rife{200415.01}, it only remains to show that there exists 
$\varepsilon_0 = \varepsilon_0 (\delta_0) > 0$ such that $0 < \| u_{\rho} - u_{\lambda} \|_{H^{1}_{0} (\Omega)} <  \delta_0$ 
whenever $0 < \| \rho - \rho_{0,\lm} \|_{\Phi}  <  \varepsilon_0$. We first prove a stronger property which will be needed in the proof of \rife{200415.1} as well. 
By using the Green's representation formula and the H\"{o}lder's inequality \rife{HoldOrl} we see that 
\beq\label{200415.5}
\|u_{\rho} - u_\lambda \|_{L^{\infty}(\om)} \leq \sup\limits_{x\in\om} \ino G(x,y)|\rho (y)-\rho_{0,\lm} (y)| \,dy \leq
C_0 \| \rho- \rho_{0,\lm} \|_{\Phi}
\eeq
where $C_0$ is the constant found in \rife{200415.3}.
Next observe that
$$\ino |\grad (u_{\rho}- u_\lambda)|^2 = \ino (|\grad u_{\rho}|^2 - 2(\grad u_{\rho},\grad u_\lambda) + |\grad u_\lambda|^2)=  
\ino \rho (u_{\rho}-u_\lambda) + \ino u_\lambda (\rho_{0,\lm}-\rho)\leq $$
$$C_0 \|\rho - \rho_{0,\lm} \|_{\Phi} + \| u_\lambda \|_{\Psi} \|\rho - \rho_{0,\lm} \|_{\Phi} \leq C \|\rho - \rho_{0,\lm} \|_{\Phi},$$
where we used \rife{130613.1}, \rife{HoldOrl}, \rife{200415.5} and the fact that obviously $\| u_\lambda \|_{\Psi}$ is bounded. We conclude that there exists $C>1$ such that
\beq\label{norm}
\|u_{\rho}-u_\lm\|_{H^{1}_0(\om)} \leq C \|\rho - \rho_{0,\lm} \|_{\Phi},
\eeq
and in particular that it is always possible to fix $\varepsilon_0:=\frac{\dt_0}{2C}>0$ so that $0 \leq \|u_\rho- u_\lambda\|_{H^{1}_0(\om)}< \dt_0$ 
whenever $0<\|\rho - \rho_{0,\lm} \|_{\Phi}<\varepsilon_0$.
Hence  \rife{200415.01} follows whenever we can prove that if $\rho$ has been chosen in this way, and therefore does not coincide with $\rho_{0,\lm}$, 
then the unique $u_{\rho}$ determined through \eqref{conc}
does not coincide with $u_\lambda$. However this is easily verified since if this was not the case we would find
$$0=-\Delta(u_{\rho}- u_\lambda) = (\rho -\rho_{0,\lm})\quad \mbox{in}\;\om$$
which is in contradiction with the fact that $\rho$ does not coincide with $\rho_{0,\lm}$. At this point \rife{030415.1} shows that 
$$
\mathcal{F} (\rho) - \mathcal{F} (\rho_\lambda) \geq J_\lambda (u_{\rho}) - J_\lambda (u_\lambda)>0,
$$
whenever $0<\|\rho - \rho_{0,\lm} \|_{\Phi}<\varepsilon_0$ as claimed.\\ 

Concerning \rife{200415.1} we first observe that, in view of \rife{l1l2}, we have $L^2(\om)\subset L_{\Phi}(\om)$. The linear operator $T:L^2(\om)\mapsto H^1_0(\om)$ 
which maps $\rho\in L^2(\om)$ to the unique weak solution $u_\rho\in H^1_0(\om)$ of \rife{conc} is a continuous bijection, whence there exists $C_3>0$ such that 
$$
\|\rho\|_{L^2(\om)}\leq C_3 \|u_\rho\|_{H^1_0(\om)},\quad \fo\,\rho\in L^2(\om),
$$
and then, since $\rho_{0,\lm}\in L^{2}(\om)$, for any $\rho\in L^{2}(\om)\cap \mathcal{P}$  we find
\beq\label{210415.1}
\|\rho-\rho_{0,\lm}\|_{L^2(\om)}\leq C_3\|u_\rho-u_\lm\|_{H^1_0(\om)}.
\eeq

On the other side \rife{l1l2} implies 

\beq\label{210415.2}
\ino \Phi(|\rho-\rho_{0,\lm}|)\leq \ino |\rho-\rho_{0,\lm}|^2,
\eeq

and then we have
\ble For any $\varepsilon_1<\varepsilon_0$ there exists $\dt_{2}>0$ such that 
\beq\label{210415.3}
\|\rho-\rho_{0,\lm}\|_{L^2(\om)} \geq \dt_2,
\eeq
for any $\rho\in L^{2}(\om)\cap \mathcal{P}$ such that $\|\rho - \rho_{0,\lm} \|_{\Phi}=\varepsilon_1$.
\ele
\proof
If the claim were false we could find a sequence $\rho_n$ such that $\|\rho_n - \rho_{0,\lm} \|_{\Phi}=\varepsilon_1$ and 
$\|\rho_n - \rho_{0,\lm}\|_{L^2(\om)} \leq \frac1n$, $\ainf$. Then \rife{210415.2} implies that $\ino \Phi(|\rho_n-\rho_{0,\lm}|)\to 0$, $\ainf$ which in view of 
\rife{170415.1} is the same as \rife{190415.1}, that is $\|\rho_n - \rho_{0,\lm} \|_{\Phi}\to 0$, $\ainf$. This is the desired contradiction to $\|\rho_n - \rho_{0,\lm} \|_{\Phi}=\varepsilon_1$
\fineproof

\bigskip

At this point, by using \rife{norm},  \rife{210415.1} and \rife{210415.3} we conclude that for any $\varepsilon_1<\varepsilon_0$ we have
$$
\frac{\dt_2}{C_3}\leq \|u_\rho-u_\lm\|_{H^1_0(\om)}\leq C\varepsilon_1,
$$
for any $\rho\in L^{2}(\om)\cap \mathcal{P}$ such that $\|\rho - \rho_{0,\lm} \|_{\Phi}=\varepsilon_1$. 
Then, in particular, we can  choose $\varepsilon_1<\frac{\dt_0}{2C}$ and a smaller $\dt_2$ if needed which satisfies $\dt_{0,1}:=\frac{\dt_2}{C_3}<C\varepsilon_1$ 
to conclude that 
\beq\label{120615.1}
\dt_{0,1}\leq \|u_\rho-u_\lm\|_{H^1_0(\om)}<\dt_0,
\eeq
for any $\rho\in L^{2}(\om)\cap \mathcal{P}$ such that $\|\rho - \rho_{0,\lm} \|_{\Phi}=\varepsilon_1$. At this point we can prove the following\\

{\bf Claim:}
{\it 
There exists $0<\eta_{1}<\eta_0:=\dt_0$ such that
$$
\eta_{1}\leq \|u_\rho-u_\lm\|_{L^2(\om)}<\eta_0,
$$
for any $\rho\in L^{2}(\om)\cap \mathcal{P}$ such that $\|\rho - \rho_{0,\lm} \|_{\Phi}=\varepsilon_1$.}\\

{\it Proof of Claim:}\\
Clearly \rife{120615.1} holds and then in particular we see that $ \|u_\rho-u_\lm\|_{L^2(\om)}<\dt_0$. Concerning the left hand side inequality we argue 
by contradiction. If the claim were false we could find a sequence $\rho_n$ such that, setting $u_n=u_{\rho_n}$, then, in view of \rife{210415.3}, we would find 
\beq\label{120615.4}
\|u_n-u_\lm\|_{H^1_0(\om)} \geq \dt_{0,1}\quad \mbox{and}\quad\|u_n-u_\lm\|_{L^2(\om)}\leq \frac1n,\,\fo\,n\in\mathbb{N}.
\eeq

On the other side we have $-\Delta (u_n-u_\lm)=(\rho_n-\rho_\lm)$ and then, multiplying by $u_n-u_\lm$ and integrating by parts, we find
$$
\|\nabla (u_n-u_\lm)\|^2_{L^2(\om)} =\ino (\rho_n-\rho_{0,\lm})(u_n-u_\lm)\leq \|\rho_n-\rho_{0,\lm}\|_{L^2(\om)}\|u_n-u_\lm\|_{L^2(\om)}\leq 
$$
$$
C_3\|u_n-u_\lm\|_{H^1_0(\om)}\frac1n\leq C_3C\|\nabla (u_n-u_\lm)\|_{L^2(\om)} \frac1n,\,\fo\,n\in\mathbb{N},
$$
where we used \rife{210415.1} and the Sobolev's inequality to conclude that $\|u_n-u_\lm\|_{H^1_0(\om)}\leq C \|\nabla (u_n-u_\lm)\|_{L^2(\om)}$, for some uniform constant $C>0$.
Therefore we would have $\|\nabla (u_n-u_\lm)\|^2_{L^2(\om)} \to 0$ and then in particular $\|u_n-u_\lm\|_{H^1_0(\om)}\to 0$, as $n\to +\infty$, which 
is a contradiction to \rife{120615.4}.
\fineproof

\bigskip

By using the statement of the Claim, then \rife{030415.1} and \rife{220415.1} imply 
that \rife{200415.1} holds. \fineproof
\bigskip
\bigskip

\section{A global existence result for \eqref{system}}\label{s2}

\bigskip

With the aid of property ${\bf (H)_{\lm}}$ in Proposition \ref{min} a global existence result follows by a standard Lyapunov 
stability argument to be applied to a set of suitably chosen initial data. 
\bpr\label{pr1} Suppose that  ${\bf (H)_{\lm}}$ in Proposition \ref{min} holds for some $\lm>0$. Then there exists $\varepsilon_\lm>0$ such that if $\rho_0$ 
in \eqref{system} 
is any smooth and non negative density such that $\| \rho_0 - \rho_{0,\lm} \|_{\Phi} \leq \varepsilon_\lm$, then $\frac{\lm}{2}\leq \ino \rho_0\leq 2\lm$
and the corresponding solution $(\rho_\lm(\cdot,t),u_\lm(\cdot,t))$ is global and uniformly bounded. 
\epr
\proof 
Since $\mathcal{F}$ is continuous, we can choose $0 < \varepsilon_2 < \varepsilon_1$ such that 
$$
\mathcal{F} (\rho) \leq \mathcal{F} (\rho_{0,\lm}) + \frac{d_1}{2},
$$ 
for any $\| \rho - \rho_{0,\lm} \|_{\Phi} \leq \varepsilon_2$. By taking a smaller value of $\varepsilon_2$ we may assume that 
$\frac{\lm}{2}\leq \ino \rho\leq 2\lm$ whenever $\| \rho - \rho_{0,\lm} \|_{\Phi} \leq \varepsilon_2$. 
At this point let us choose $\rho_0$ in \eqref{system} 
to be any smooth and non negative density such that 
$\| \rho_0 - \rho_{0,\lm} \|_{\Phi} \leq \varepsilon_2$ and let $(\rho_\lm,u_\lm)$ denote the corresponding unique solution in $[0,T]$, 
for some $T>0$.
Then, we claim that for any $t\in[0,T]$ it holds
\beq\label{main}
\mathcal{F} (\rho_{0,\lm}) \leq \mathcal{F} (\rho (\cdot ,t)) < \mathcal{F} (\rho_{0,\lm}) + d_1.
\eeq
In fact, on one side the right hand inequality is always satisfied since \eqref{dissip} implies that
$$
\mathcal{F} (\rho (\cdot ,t)) \leq \mathcal{F} (\rho_{0}) \leq \mathcal{F} (\rho_{0,\lm}) + \frac{d_1}{2}.
$$
On the other side we also have $\rho (\cdot,t)\in \mathcal{P}_{\lm}$ in view of the mass conservation and the fact that $\rho(x,t)$ is 
classical whence it satisfies \rife{global} whenever it exists 
(see \cite{suzC} Theorem 3.1). So, if for some $t_*>0$ we would find that $\rho_{*}:=\rho(x, t_*)$ satisfies  
$\mathcal{F} (\rho_{*}) < \mathcal{F} (\rho_{0,\lm})$, 
then, in view of ${\bf (H)_{\lm}}$, necessarily $\| \rho_{*} - \rho_{0,\lm} \|_{\Phi} \geq \varepsilon_0$. Therefore in particular, by the continuity of the norm, 
there existed $t_1<t_*$ such that $\rho_{1}:=\rho(x, t_1)$ 
satisfied $\| \rho_{1} - \rho_{0,\lm} \|_{\Phi} =\varepsilon_1$. But then \rife{200415.1} 
implies that $\mathcal{F} (\rho_{1}) \geq \mathcal{F} (\rho_{0,\lm})+d_1$, since $\rho_1$ being classical surely belongs to $L^2(\om)$. This is 
the desired contradiction and thus we have proved that \rife{main} holds. In particular it is easy to see that \rife{main} implies that there exist $C>0$, 
depending only 
on $\lm$ and $d_1$, such that
\beq\label{main1}
\ino \rho_\lm(\cdot,t)(\log{(\rho_\lm(\cdot,t))}-1)\leq C,\qquad \ino |\nabla u_\lm(\cdot, t)|^2\leq C,
\eeq
for any $t\in [0,T]$. At this point we can follow step by step the argument in Theorem 3 in \cite{suz} to conclude that 
$$
\sup\limits_{t\in [0,T]}\sup\limits_{x\in\ov{\om}}\rho_\lm(x,t)\leq \tilde{C},
$$
for some uniform constant $\tilde{C}>0$. Then well known arguments imply that the solution is global and the desired conclusion follows
by choosing $\varepsilon_\lm=\varepsilon_2$.
\fineproof

\bigskip
\bigskip

\section{The Proof of Theorem \ref{kstype}}\label{theoproof}
In this section we prove Theorem \ref{kstype}.

{\bf The Proof of Theorem \ref{kstype}}\\
We discuss the proof of part {\rm (a)}. The proof of part {\rm (b)} can be worked exactly with the same argument with minor changes.\\
By Theorem \ref{subsup} we have a strict local minimizer $u_\lm$ of $J_\lambda$ for any $\lm\leq \lm_{\a,c_{\sscp D}}$, where 
$\underline{\lm}_{\a,c_{\sscp D}}< \lm_{\a,c_{\sscp D}}<\overline{\lm}_\a$ and 
$\underline{\lm}_{\a,c_{\sscp D}}\simeq\frac{4\pi c_{\sscp D}}{(8-c_{\sscp D})\a}$, $\ov{\lm}_{\a}\simeq
\frac{2 \pi}{3\a}$ as $\a\to 0^+$. It follows from Proposition \ref{min} that  if $\rho_{0,\lm}$ is defined as in \rife{030415.2}, then 
$\rho_{0,\lm}\in \mathcal{P}_\lm$ is a strict local minimizer of $\mathcal{F}$.  Actually we have a stronger result since ${\bf (H)_{\lm}}$ in 
Proposition \ref{min} holds. Therefore we can apply Proposition \ref{pr1} to conclude that 
there exists $\varepsilon_\lm>0$ such that if $\rho_0$ in \eqref{system} 
is any smooth and non negative density such that $\| \rho_0 - \rho_{0,\lm} \|_{\Phi} \leq \varepsilon_\lm$, then $\frac{\lm}{2}\leq \ino \rho_0\leq 2\lm$
and the corresponding solution $(\rho_\lm,u_\lm)$ is global and uniformly bounded. Let $2 m_\lm=\min\limits_{\ov{\om}} \rho_{0,\lm}$. Clearly $m_\lm>0$ 
and we define $f_\lm$ to be any smooth function in $\ov{\om}$ which satisfies $|f_\lm|\leq m_\lm$ and $\ino f_\lm=0$. 
Then we can choose $0<\sg<\frac{1}{2}$ depending 
on $\varepsilon_\lm$ such that 
$\rho_0=\rho_{0,\lm}+\sg f_{\lm}$ satisfies $0<\| \rho_0 - \rho_{0,\lm} \|_{\Phi} \leq \varepsilon_\lm$. Clearly $\ino \rho_{0}=\lm$ and then 
in particular $\rho_0\in \mathcal{P}_\lambda$. Therefore, 
for any $\lm\leq \lm_{\a,c_D}$ we have found initial data $\rho_0$ such that the solution $(\rho_\lm,u_\lm)$ of $P(\lm,\om)$ is global and uniformly bounded 
as claimed.
\fineproof

\end{document}